\documentclass[12pt,a4paper]{amsart}

\newtheorem{theo+}              {Theorem}           [section]
\newtheorem{prop+}  [theo+]     {Proposition}
\newtheorem{coro+}  [theo+]     {Corollary}
\newtheorem{lemm+}  [theo+]     {Lemma}
\newtheorem{exam+}  [theo+]     {Example}
\newtheorem{rema+}  [theo+]     {Remark}
\newtheorem{defi+}  [theo+]     {Definition}
\def \r{\mbox{${\mathbb R}$}}

\newenvironment{theorem}{\begin{theo+}}{\end{theo+}}

\usepackage{amsthm}
\theoremstyle{plain} \theoremstyle{remark}

\def\E{/\kern-1.0em \equiv }

\title{Some recent progress of biharmonic submanifolds}
\author{Ye-Lin Ou $^{*}$}
\thanks{$^{*}$ Research supported by
Faculty Development Program of Texas A $\&$ M University-Commerce, 2015.}
\address{Department of
Mathematics,\newline\indent Texas A $\&$ M University-Commerce,
\newline\indent Commerce, TX 75429,\newline\indent USA.\newline\indent
E-mail:yelin.ou@tamuc.edu }
\begin{document}
\title[Some recent progress of biharmonic submanifolds]{Some recent progress of biharmonic submanifolds}
\date{11/28/2015}
\maketitle
\section*{Abstract}
\begin{quote} In this note, we give a brief survey on some recent developments of biharmonic submanifolds. After reviewing some recent progress on Chen's biharmonic conjecture, the Generalized Chen's conjecture on biharmonic submanifolds of non-positively curved manifolds, and some classifications of biharmonic submanifolds of spheres, we give a short list of some open problems in this topic.
{\footnotesize } 
\end{quote}
\section{Biharmonic submanifolds and their equations} 
{\bf 1.1 Biharmonic submanifolds}\\

The study of biharmonic submanifolds began with the independent works of B. Y. Chen \cite{Ch1} and G. Y. Jiang \cite{Ji1}, \cite{Ji2}, \cite{Ji3} in the middle of 1980s. In his program to study the finite type submanifolds of Euclidean space, Chen defined a biharmonic submanifold of a Euclidean space to be a submanifold with harmonic mean curvature vector field, i.e., an isometric immersion $\phi: M^m\longrightarrow \mathbb{R}^n$ that satisfies $\Delta H=0$, where $H$ is the mean curvature vector field and $\Delta$ denotes the Laplacian on the submanifold $M$. In the meantime, G. Y. Jiang, in his effort to study biharmonic maps between Riemannian manifolds as a part of a program to understand the geometry of $k$-polyharmonic maps proposed by Eells and Lemaire in \cite{EL2}, began to study biharmonic submanifolds of Riemannian manifolds as biharmonic isometric immersions. The paper \cite{CMO1} of Caddeo-Montaldo-Oniciuc seems to be the first one to use the term of ``biharmonic submanifolds of Riemannian manifolds" .\\

A {\em harmonic map} is a map $\varphi:(M, g)\longrightarrow (N, h)$ between Riemannian manifolds that is a critical point of the energy functional
\begin{equation}
E\left(\varphi \right)= \frac{1}{2} {\int}_{M} \left|d \varphi
\right|^{2}dv_{g},
\end{equation}
where $\frac{1}{2}\left|{\rm
d}\varphi\right|^{2}=\frac{1}{2}g^{ij}\varphi^{\alpha}_{i}
\varphi^{\beta}_{j}h_{\alpha \beta}$ is the energy density of $\varphi$. Harmonic map equation is a system of 2nd order PDEs given by the Euler-Lagrange equation of the energy functional:
\begin{align}\notag
\Delta \varphi^{\sigma}+g(\nabla \varphi^{\alpha}, \nabla\varphi^{\beta})\bar{\Gamma}^{\sigma}_{\alpha
\beta}=0, \;\;\;\; \sigma =1, 2, \cdots, n.
\end{align}
The tension field of the map $\varphi$ is defined to be 
\begin{align}\notag
\tau(\varphi):={\rm Trace}_g\nabla d\varphi
=g^{ij}(\varphi_{ij}^{\sigma}-\Gamma_{ij}^{k}\varphi_{k}^{\sigma}+
\bar{\Gamma}^{\sigma}_{\alpha
\beta}\varphi^{\alpha}_{i}\varphi^{\beta}_{j})\frac{\partial}{\partial
y^{\sigma}}.
\end{align}
So, harmonic maps are those maps whose tension fields vanish identically.\\

It is well known that harmonic maps include many important classes of maps in mathematics as special cases, including harmonic functions, geodesics, and minimal isometric immersions (i.e., minimal submanifolds). In fact, for an isometric immersion $\varphi:(M^m, g)\longrightarrow (N^n, h)$, its tension field is $\tau(\varphi)=m H$, where $H$ is the mean curvature vector field of the submanifold $\varphi(M) \subseteq (N^n, h)$. So an isometric immersion is a harmonic map if and only if $\varphi(M) \subseteq (N^n, h)\;\;{\rm is\; a\;minimal \;submanifold}$ (cf. \cite{ES} ).\\

A {\em biharmonic map} is a map $\varphi:(M, g)\longrightarrow (N, h)$ between Riemannian manifolds that is a critical point of the bi-energy functional
\begin{equation}
E_2\left(\varphi \right)= \frac{1}{2} {\int}_{M} \left|\tau( \varphi)
\right|^{2}dv_{g},
\end{equation}
where $\tau(\varphi)=
{\rm Trace}_g\nabla d\varphi$ is the tension field of $\varphi$.
Biharmonic map equation is a system of a 4th order PDEs given by the Euler-Lagrange equation of the bienergy functional (See \cite{Ji2}):
\begin{equation}\notag
\tau_{2}(\varphi):={\rm
Tr}_{g}(\nabla^{\varphi}\nabla^{\varphi}-\nabla^{\varphi}_{\nabla^{M}})\tau(\varphi)
- {\rm Tr}_{g} R^{N}({\rm d}\varphi, \tau(\varphi)){\rm d}\varphi
=0,
\end{equation}
where
$R^{N}(X,Y)Z=
[\nabla^{N}_{X},\nabla^{N}_{Y}]Z-\nabla^{N}_{[X,Y]}Z$ is the Riemann curvature operator of $(N, h)$.\\

{\em A biharmonic submanifold} is a submanifold whose defining isometric immersion is a biharmonic map. It is clear from the definitions that a harmonic map is always a biharmonic map and hence a minimal submanifold is always a biharmonic submanifold. It is customary to call a biharmonic map which is not harmonic {\em a proper biharmonic map} and a biharmonic submanifold that is not minimal {\em a proper biharmonic submanifold}.\\

{\bf 1.2 The equations of biharmonic submanifolds}\\

The equations for a biharmonic submanifold of a Riemannian manifold in local coordinates were first derived by Jiang in \cite{Ji1} as
\begin{equation}
\begin{cases}
\sum_{j,k=1}^m\{ b^{\alpha}_{jjkk}-\sum_{\beta=m+1}^nb^{\beta}_{jj}b^{\beta}_{lk}b^{\alpha}_{lk}+b^{\beta}_{jj}R^{N}_{\beta k \alpha k}\}=0,\\
\sum_{j,k=1}^m \sum_{\beta=m+1}^n\{ 2 b^{\alpha}_{jjk}b^{\beta}_{ik}+b^{\beta}_{jj}b^{\beta}_{ikk}-b^{\beta}_{jj}R^{N}_{\beta k i k}\}=0,\\\alpha=m+1, m+2, \cdots, n, i=1, 2, \cdots, m,
\end{cases}
\end{equation}
where $R^{N}_{abcd}$, $b^{\alpha}_{ij}$, $b^{\alpha}_{ijk}$, $b^{\alpha}_{ijk}$ denote the curvature tensor of the ambient space $N$, the second fundamental form of the submanifold, the first and the second covariant derivatives of $b^{\alpha}_{ij}$ respectively. Our notation convention for the curvature is such that the Weingarten formula reads $b^{\alpha}_{ijk}-b^{\alpha}_{ikj}=R^N_{\alpha i\,kj}$.\\

Invariant form of the equations for biharmonic submanifold of Euclidean and Semi-Euclidean spaces were derived and used by Chen in \cite{Ch0}, \cite{Ch1}, and that of space forms of nonzero constant curvatures were derived and used by Oniciuc in \cite{On1} and Caddeo-Montaldo-Oniciuc in \cite{CMO2}.\\

The invariant form of the equations for biharmonic hypersurfaces $\varphi:M^{m}\longrightarrow N^{m+1}$ in a generic Riemannian manifold was derived by the author in \cite{Ou4}:
\begin{equation}\label{BHEq}
\begin{cases}
\Delta H-H |A|^{2}+H{\rm
Ric}^N(\xi,\xi)=0,\\
2A\,({\rm grad}\,H) +\frac{m}{2} {\rm grad}\, H^2
-2\, H \,({\rm Ric}^N\,(\xi))^{\top}=0,
\end{cases}
\end{equation}
where ${\rm Ric}^N : T_qN\longrightarrow T_qN$ denotes the Ricci
operator of the ambient space defined by $\langle {\rm Ric}^N\,
(Z), W\rangle={\rm Ric}^N (Z, W)$, $H$ is the mean curvature function, and $A$ is the shape operator of the hypersurface.\\

Finally, the invariant form of biharmonic submanifold equations for a general Riemannian manifold was given in \cite{BMO7}:
\begin{equation}\label{BHEq}
\begin{cases}
\Delta^{\perp} H+ {\rm tr}\, B(\cdot, A_H\cdot)+ {\rm tr}\,
({\rm R}^N(d\varphi(\cdot),H)d\varphi(\cdot))^{\perp}=0,\\
2\,{\rm tr} A_{\nabla^{\perp}(\cdot) H}(\cdot)+\frac{m}{2} {\rm grad}\, |H|^2
+2\, {\rm tr}\,
({\rm R}^N(d\varphi(\cdot),H)d\varphi(\cdot))^{\top}=0,
\end{cases}
\end{equation}
where $H$, $\Delta^{\perp}$ and $\nabla^{\perp}$ denote the mean curvature vector field, the Laplacian and the connection of the normal bundle of the submanifold.\\

{\bf 1.3 Some examples of proper biharmonic submanifolds}\\

In general, it is not easy to find examples of proper biharmonic submanifolds due to the nature of 4th order and nonlinearity of the biharmonic submanifold equations. The following list gives some examples of proper biharmonic submanifolds.
\begin{itemize}
\item Biharmonic hypersurfaces:
\begin{itemize}
\item[(i)] The generalized Clifford tori \cite{Ji2}: $S^p(\frac{1}{\sqrt{2}})\times S^q(\frac{1}{\sqrt{2}})\hookrightarrow S^{n+1}$ with $p\ne q, \;\;p+q=n$;
\item[(ii)] Small spheres \cite{CMO1}: $S^n(\frac{1}{\sqrt{2}}) \hookrightarrow S^{n+1}$; ;
\item[(iii)] Hypercylinders (\cite{Ou4}, \cite{Ou5}): $S^{n-1}(\frac{1}{\sqrt{2}})\times \mathbb{R}\hookrightarrow S^{n}\times \mathbb{R}$;
\item[(iv)] The hyperplanes $x_{n+1}=k$ in the conformally flat spaces(\cite{Ou4}):\\ $\left(\r_{+}^{n+1}, (\frac{x_{n+1}+C}{D})^{2}\sum_{i=1}^{n+1}{\rm d}{x_i}^{2}\right)$.
\item[(v)] Some product of spheres (cf. \cite{BMO}): $S^{n_1}(r_1)\times\ldots \times S^{n_q}(r_q)\times S^{n_{q+1}}(r_{q+1})\times\ldots \times
S^{n_k}(r_k)$ in $S^{m+k-1}$, where $m=\sum_{i=1}^{k} n_i, r_i=\sqrt{n_i}/\sqrt{2(n_1+\cdots +n_q)}$ for $i=1, 2, \cdots, q$ and $r_j=\sqrt{n_j}/\sqrt{2(n_{q+1}+\cdots +n_k)}$ for $j=q+1, 2, \cdots, k$.
\item[(vi)] Some biharmonic real
hypersurfaces in $CP^n$ and some biharmonic tori:
$T^{n+1}=S^{1}(r_1)\times S^{1}(r_2)\times\ldots \times
S^{1}(r_{n+1})$ in $S^{2n+1}$ (See \cite{Zh} for details).
\end{itemize}
\item Biharmonic submanifolds via constructions of precompositions or post compositions: 
\begin{itemize}
\item[(i)] Caddeo-Montaldo-Oniciuc \cite{CMO2}: The submanifold $M^m\hookrightarrow S^{n+1}$ obtained as the composition $M^m\hookrightarrow S^n(\frac{1}{\sqrt{2}}) \hookrightarrow S^{n+1}$, where $M^m\hookrightarrow S^n(\frac{1}{\sqrt{2}})$ is any minimal submanifold; 
\item[(ii)] Caddeo-Montaldo-Oniciuc \cite{CMO2}: The submanifold $M^k\times N^l \hookrightarrow S^{n+1}$ obtained as the composition $M^k\times N^l \hookrightarrow S^p(\frac{1}{\sqrt{2}})\times S^q(\frac{1}{\sqrt{2}})\hookrightarrow S^{n+1}$ with $p\ne q, \;\;p+q=n$, where $M^k\hookrightarrow S^p(\frac{1}{\sqrt{2}})$ and $ N^l \hookrightarrow S^q(\frac{1}{\sqrt{2}})$ are two minimal submanifolds.
\item[(iii)] The submanifold $M^k\hookrightarrow S^n$ obtained as the composition $M^k\hookrightarrow S^m\hookrightarrow S^n$, where $M^k\hookrightarrow S^m$ is a proper biharmonic submanifold and $ S^m\hookrightarrow S^n$ is a totally geodesic immersion (See \cite{Ou2} and \cite{Su}).
\end{itemize}
\item A biharmonic Legendre surface in $S^5$ \cite{Sa}: $\phi: M^2\longrightarrow S^5$ with $\phi(u, v)=\\\frac{1}{\sqrt{2}}( \cos u, \sin u \sin (\sqrt{2}v), -\sin u \cos (\sqrt{2}v), \sin u, \cos u \sin (\sqrt{2}v), -\cos u \cos (\sqrt{2}v))$.
\end{itemize}
\section{Biharmonic submanifolds of nonpositively curved spaces}
{\bf 2.1 Chen's conjecture on biharmonic submanifolds of Euclidean spaces}\\

One of the fundamental problems in the study of biharmonic submanifolds is to classify such submanifolds in a model space. So far, most of the work done has been focused on classification of biharmonic submanifolds of space forms. For the case of Euclidean ambient space, all the work done has been trying to solve the following\\
{\em B. Y. Chen's conjecture on biharmonic submanifolds} (See e.g., \cite{Ch1}, \cite{Ch}, \cite{Ch2}, \cite{Ch3}, \cite{Ch4}, \cite{Ch5}, \cite{CMO1}, \cite{BMO2}, \cite{KU}, \cite{NU1}, \cite{NU2}, \cite{Ou5}, \cite{Wh} and the references therein): Any biharmonic submanifold in a Euclidean space is a minimal one.\\

$\clubsuit$ The conjecture is still open although it has been verified to be true in the following special cases:
\begin{itemize}
\item[(I)] Biharmonic curves in $\mathbb{R}^m$ (Jiang 1986 \cite{Ji3} and independently Dimitri$\acute{\rm c}$ 1989 \cite{Di1}, 1992 \cite{Di2} ).
\item[(II)] Biharmonic hypersurfaces $M^m\hookrightarrow \mathbb{R}^{m+1}$ with
\begin{itemize}
\item[$\bullet$] $m=2$ (Chen 1985 \cite{Ch1}, 2011 \cite{Ch2} and independently Jiang 1986 \cite{Ji3}).
\item[$\bullet$] $m=3$ (Hasanis-Vlachos 1995 \cite{HV} and Defever 1998 \cite{De} for a different proof).
\item[$\bullet$] parallel mean curvature ( Jiang 1986 \cite{Ji3}).
\item[$\bullet$] completeness and Ricci curvature bounded from below (Jiang 1986 \cite{Ji3}).
\item[$\bullet$] at most 2 distinct principle curvatures (Dimitri$\acute{\rm c}$ 1992 \cite{Di2}).
\item[$\bullet$] $3$ or $4$ special distinct principal curvatures: (i) $k_1, k_2, k_1+k_2$, or (ii) $k_1, k_2, k_3, k_1+k_2+k_3$ (Chen-Munteanu 2013 \cite{CM}).
\item[$\bullet$] at most $3$ distinct principal curvatures: $k_1, k_2, k_3$ (Fu 2014 \cite{Fu1}, 2015 \cite{Fu2}).
\item[$\bullet$] completeness and finite total mean curvature (Nakauchi-Urakawa 2011 \cite{NU1})
\item[$\bullet$] oriented image contained in a non-degenerated open cone or being cylindrically bounded (Alias-Martinez-Rigoli 2013 \cite{AGR}).
\item[$\bullet$] weak convexity (Luo 2014 \cite{Lu1}).
\item[$\bullet$] $m$ distinct principal curvatures and $g(\nabla_{e_i}e_k, e_j)\ne 0$ for all distinct triples of principal curvature vectors $\{ e_i, e_j, e_k\}$ in the kernel of $dH$ (Koiso-Urakawa 2014 \cite{KU}).
\item[$\bullet$] $G$-invariant property with cohomogeneity 2 group of isometry of $\r^{m+1}$ for $m\ge 2$ (See the most recent work of Montaldo-Oniciu-Ratto 2015, \cite{MOR}).
\end{itemize}
\end{itemize}
\begin{itemize}
\item[(III)] Biharmonic submanifolds $M^m\hookrightarrow \mathbb{R}^{m+p}$ which 
\begin{itemize}
\item[$\bullet$] are of finite type ( Dimitri$\acute{\rm c}$ 1989 \cite{Di1}, 1992 \cite{Di2}).
\item[$\bullet$] are pseudo-umbilical with $m\ne 4$ (Dimitri$\acute{\rm c}$ 1992 \cite{Di2}).
\item[$\bullet$] are spherical (B. -Y. Chen 1996 \cite{Ch}).
\item[$\bullet$] are properly immersed (Akutagawa-Maeta 2013 \cite{AM}).
\item[$\bullet$] are complete and have finite total mean curvature (Nakauchi-Urakawa 2013 \cite{NU2})
\item[$\bullet$] satisfy the decay condition at infinity (Wheeler 2013 \cite{Wh}).
\item[$\bullet$] are complete and have at most polynomial growth (Luo 2015 \cite{Lu2}).
\end{itemize}
\end{itemize}
{\bf 2.2 Biharmonic graphs and similarity between Bernstein and Chen's Conjectures}\\

It is well known that any hypersurface is locally the graph of a function. As minimal graphs have play an important role in the study of minimal (hyper)surfaces, the author initiated the study of biharmonic graphs in \cite{Ou5}. 
\begin{theorem} \cite{Ou5} 
$( I ) $ A function $f:\mathbb{R}^m\supseteq D\longrightarrow \mathbb{R}$ has
biharmonic graph \\$\{(x, f(x)): x\in D\} \subseteq \mathbb{R}^{m+1}$ if and
only if  $f$ is a solution of
\begin{equation}\label{Bg1}
\begin{cases}
\Delta (\Delta f)=0,\\(\Delta f_k)\Delta f+2g(\nabla f_k,\nabla \Delta
f)=0, \;\; k=1,2,\ldots, m,
\end{cases}
\end{equation}
where the Laplacian $\Delta$ and the gradient $\nabla$ are taken with respect to the induced metric $g_{ij}=\delta_{ij}+f_if_j$.\\
$(II)$ Chen's Conjecture for biharmonic hypersurfaces is equivalent to stating that any solution of Equation (\ref{Bg1}) is a harmonic function or, equivalently, the only solution of (\ref{Bg1}) is the ``trivial" one---the one that satisfies $\Delta f=0$.
\end{theorem}
Recall that the important {\bf Bernstein Conjecture} in the history of mathematics states that for $m\ge 2$, any entire solution $f: \r^m \longrightarrow \r$ of the minimal graph equation
\begin{equation}\notag
\Delta f=0 \;\;\;\;\;\;\Longleftrightarrow \;\;\;\;\;\;
\sum_{i,j=1}^{m}\left(\delta_{ij}-\frac{f_{i}f_{j}}{1+\left|\nabla
f\mathbb{R}\right|^{2}}\mathbb{R}\right)f_{ij}=0,
\end{equation}
is the ``trivial" one---the one with $f_{ij}=0$, i.e., an affine function. In this sense, we see a similarity between Bernstein and Chen's conjectures.\\

It is also well known that Bernstein conjecture
\begin{itemize}
\item is true for \begin{itemize}
\item[] $m=2$ by Bernstein's work (1915);
\item[] $m=3$ by De Giorgi's work (1965);
\item[] $m=4$ by Almgren's work (1966);
\item[] $m\leq 7$ by J. Simons's work (1968), but 
\end{itemize}
\item is false for $m>7$ by the joint work of Bombieri-De Giorgi-Giusti (1969);
\end{itemize}
Now for Chen's Conjecture for biharmonic graphs of $f: \r^m\supset \Omega \longrightarrow \r$, we know that
\begin{itemize}
\item It is true for $m=2$ by the work of Chen (1985) and Jiang (1986) independently;
\item It is true for $m=3$ by the work of Hasanis-Vlachos (1995).
\end{itemize}
Could it be possible that Chen's conjecture fails to be true for some $m \ge 4$ similar to the situation for the Bernstein conjecture?\\

{\bf 2.3 Biharmonic Submanifolds of nonpositively curved spaces}\\

An interesting problem related to Chen's conjecture of biharmonic submanifolds is the following generalized Chen's conjecture for biharmonic submanifolds proposed by Caddeo, Montaldo, and Oniciuc in 2001.\\

{\em The Generalized Chen's conjecture of biharmonic submanifolds} ( See \cite{CMO1}, \cite{BMO2}, also \cite{MO}, \cite{Ch2}, \cite{Ch3}, \cite{Ch5}, \cite{IIU2}, \cite{Ou4}, and \cite{OT} etc): Any biharmonic submanifolds of nonpositively curved manifold is a minimal one.\\

Special cases that support the Generalized Chen's conjecture on biharmonic submanifolds include:
\begin{itemize}
\item Compact submanifolds (Jiang 1986 \cite{Ji2});
\item Biharmonic submanifold of $H^3(-1)$ (Caddeo-Montaldo-Oniciuc 2002 \cite{CMO2});
\item Biharmonic curves of $H^n(-1)$ (Caddeo-Montaldo-Oniciuc 2002 \cite{CMO2});
\item Pseudo-umbilical biharmonic submanifold $M^m\subset H^n(-1)$ with $m\ne 4$ (Caddeo-Montaldo-Oniciuc 2002 \cite{CMO2}); 
\item Biharmonic submanifold  of constant mean curvature (Oniciuc 2002 \cite{On1}); 
\item Biharmonic hypersurface of $H^n(-1)$ with at most two distinct principal curvatures (Balmus-Montaldo-Oniciuc 2008 \cite{BMO2}).
\end{itemize}

However, the generalized Chen's conjecture was proved to be false by Ou and Tang in \cite{OT}. Ou and Tang first constructed a family of counter examples of $4$-dimensional proper biharmonic hypersurfaces in a conformally flat space with negative sectional curvature \cite{OT}: \\

{\em For constants $ A>0,\; B>0,\;c $, any $t\in
 (0,\;\frac{1}{2})$, and any $(a_1, a_2, a_3,
a_4)\in S^3\left(\sqrt{\frac{2t}{1-2t}}\;\right)$, the
isometric immersion
\begin{equation}\notag
\varphi : \mathbb{R}^4\longrightarrow
\left(\mathbb{R}^5_{+}, (Az+B)^{-2t}[\sum_{i=1}^4{\rm d}{x_i}^{2}+{\rm
d}{z}^{2}]\right)
\end{equation}
with $\varphi(x_1,\ldots,
x_4)=(x_1,\ldots, x_4,\sum_{i=1}^{4}a_ix_i+c)$ is proper biharmonic
into the conformally flat space. Furthermore, for $t\in (0,1)$, the conformally flat
space $(\mathbb{R}_{+}^{5},h=(Az+B)^{-2t}(\sum_{i=1}^4{\rm
d}{x_i}^{2}+{\rm d}{z}^{2}))$ has negative sectional curvature.
}\\

By using the product maps of these family of counter examples and totally geodesic immersions, Ou and Tang then provided the following counter examples \cite{OT}:\\

 {\em For any positive integer $k$ and $n$, let $\varphi:M^4\longrightarrow (\mathbb{R}^5, h) $ be a counter example of proper biharmonic hypersurface in the conformally flat space constructed above, and $\psi: \mathbb{R}^n\longrightarrow (\mathbb{R}^{n+k-1}, h_0)$ be the standard totally geodesic embedding of Euclidean spaces. Then, the isometric immersion $\phi:M^4\times \mathbb{R}^n\longrightarrow (\mathbb{R}^5\times\mathbb{R}^{n+k-1}, h+h_0)$
with $\phi(x, y)= (\varphi(x), \psi (y))$ is a proper biharmonic submanifold of codimension $k$ in a nonpositively curved space.}

\section{Biharmonic submanifolds of spheres--some classifications}
The study of biharmonic submanifolds of spheres was initiated by Jiang in \cite{Ji1}, \cite{Ji2}. Some of Jiang's work in this area can be summarized as
\begin{itemize}
\item Derivation of the equations of biharmonic submanifolds in local coordinates in a generic Riemannian manifold and in space forms (\cite{Ji1}); \item Obtaining the following partial classifications/characterizations:
\begin{itemize}
\item a compact biharmonic hypersurface $\phi: M^m\hookrightarrow S^{m+1}$ with square norm of the second fundamental form $|B(\phi)|^2<m$ is minimal. Furthermore, if $|B(\phi)|^2=m$, then the biharmonic hypersurface has constant mean curvature (\cite{Ji1});
\item biharmonic submanifolds $\phi: M^m\hookrightarrow S^{m+p}$ with parallel mean curvature vector field and square norm of the second fundamental form $|B(\phi)|^2<m$ is minimal (\cite{Ji1});
\item a hypersurface $\phi: M^m\hookrightarrow S^{m+1}$ with parallel mean curvature is biharmonic if and only if its Gauss map is biharmonic (\cite{Ji1}); 
\item a hypersurface $\phi: M^m\hookrightarrow S^{m+1}$ with parallel mean curvature vector field and non-zero mean curvature is biharmonic if and only if $|B(\phi)|^2= m$ (\cite{Ji2}). As an application of this, Jiang \cite{Ji2} also gave the first example of proper biharmonic submanifold $S^p\times S^q\hookrightarrow S^{p+q+1}$ with $p\ne q$.
\end{itemize}
\end{itemize}

Since 2001, the study of biharmonic submanifolds of spheres has been vigorously extended by Balmus, Caddeo, Montaldo, Oniciuc, Urakawa and others. The main focus is to classify/characterize biharmonic submanifolds in a sphere. Some partial results include the following.\\
(I) Biharmonic curves in a sphere -- a complete classification has been achieved:
\begin{itemize}
\item Any proper biharmonic  curve in $S^2$ is a part of $S^1(\frac{1}{\sqrt{2}})$ (Caddeo-Montaldo-Piu \cite{CMP}).
\item Any proper biharmonic curve in $S^3$ is a part of $S^1(\frac{1}{\sqrt{2}})\subseteq S^3$, or
geodesic of $S^1(\frac{1}{\sqrt{2}})\times S^1(\frac{1}{\sqrt{2}})\subseteq S^3$
(Caddeo-Montaldo-Oniciuc \cite{CMO1}).
\item Any proper biharmonic curve in $S^n$ $( n>3)$ is a biharmonic curve in $S^3$ sitting in $S^n$ as a totally geodesic submanifolds (Caddeo-Montaldo-Oniciuc \cite{CMO2}).
\end{itemize}
(II) Biharmonic hypersurfaces -- Some partial classifications:
\begin{itemize}
\item Any proper biharmonic surface in $S^3$ is a part of $S^2(\frac{1}{\sqrt{2}})$ (Caddeo-Montaldo-Oniciuc \cite{CMO1}).
\item The only proper biharmonic compact hypersurfaces of
$S^4$ are the hypersphere $S^3(\frac{1}{\sqrt{2}})$ and the torus
$S^1(\frac{1}{\sqrt{2}})\times S^2(\frac{1}{\sqrt{2}})$ (Balmus-Montaldo-Oniciuc \cite{BMO4}).
\item A hypersurface $\phi: M^{m}\hookrightarrow S^{m+1}$ with at most two distinct principal curvatures is proper biharmonic, then $\phi(M)$ is CMC and is an open part of the hypersphere $S^{m}(\frac{1}{\sqrt{2}})$, or the generalized Clifford torus $S^p(\frac{1}{\sqrt{2}})\times S^{q}(\frac{1}{\sqrt{2}})$ with
$p\ne q, p+q=m$ (Balmus-Montaldo-Oniciuc \cite{BMO2}).
\item There does not exist compact CMC proper biharmonic hypersurfaces with $3$ distinct principal curvatures everywhere (Balmus-Montaldo-Oniciuc \cite{BMO4}).
\item If a proper biharmonic hypersurface $\phi: M^m\hookrightarrow S^{m+1}$ has constant mean curvature, then $|A|^2=m$; Furthermore for $m>2$, then $|H|\in (0, \frac{m-2}{m})\cup \{1\}$, and if $|H|=1$ or $|H|=(m-2)/m$, then $\phi(M)$ is an open part of $S^{m}(\frac{1}{\sqrt{2}})$ or $S^{m-1}(\frac{1}{\sqrt{2}})\times S^1(\frac{1}{\sqrt{2}})$ respectively (Balmus-Montaldo-Oniciuc \cite{BMO2} and \cite{BO2}).
\item An isoparametric hypersurface in
$S^{m+1}$ is proper biharmonic if and only if it is a part of either
$S^m(1/\sqrt{ 2} )$ or $S^p(1/\sqrt{ 2} )\times S^q(1/\sqrt{ 2} ),
p+q = m, p\ne q$. (Ichiyama-Ignoguchi-Urakawa \cite{IIU2}).
\item An orientable Dupin hypersurface $\phi: M^m\hookrightarrow S^{m+1}$ with at most three distinct principal curvatures is proper biharmonic if and only if it is an open part of $S^{m}(\frac{1}{\sqrt{2}})$ or $S^p(\frac{1}{\sqrt{2}})\times S^{q}(\frac{1}{\sqrt{2}})$ with
$p\ne q, p+q=m$ (Balmus-Montaldo-Oniciuc \cite{BMO6}).
\item A compact and  orientable  hypersurface $\phi: M^m\hookrightarrow S^{m+1}$ with at most three distinct principal curvatures is proper biharmonic if and only if it is either $S^{m}(\frac{1}{\sqrt{2}})$ or $S^p(\frac{1}{\sqrt{2}})\times S^{q}(\frac{1}{\sqrt{2}})$ with
$p\ne q, p+q=m$ (Fu \cite{Fu3}).
\end{itemize}

All the known examples and results of biharmonic submanifolds in spheres suggests the following\\
{\bf Conjecture} (Balmus-Montaldo-Oniciuc \cite{BMO2}, \cite{BMO3}, and \cite{BMO6}): any biharmonic submanifold in sphere has constant mean curvature; and any proper biharmonic hypersurface in $S^{n+1}$ is an
open part of the hypersphere $S^{n}(\frac{1}{\sqrt{2}})$, or the
generalized Clifford torus $S^p(\frac{1}{\sqrt{2}})\times
S^{q}(\frac{1}{\sqrt{2}})$ with $p\ne q, p+q=n$.

\section{Some open problems and final remarks}
In this final section, we would like to list the following open problems in the study of biharmonic submanifolds.\\
{\bf (1)  Chen's conjecture on biharmonic submanifold}.\\
The conjecture is still open and the author believes that the following are the interesting and critical cases to work on
\begin{itemize}
\item[C1.] Hypersurfaces $M^4\hookrightarrow \r^5$. Based on the known results, one needs only to check whether the hypersurfaces with simple principal curvatures, i.e., $4$ distinct principal curvatures is minimal or not;
\item[C2.] Pseudo-umbilical submanifold $M^4\hookrightarrow\r^{n+4}$ with $n>1$. This in the only case not covered in Dimitri$\acute{\rm c}$'s results. It is interesting to note that the number ``$4$" has mysteriously appeared in several cases associated with biharmonicity: The first counter example to the Generalized Chen's conjecture is a $4$-dimensional hypersurface (cf. \cite{OT}), biharmonic morphism has to be a $4$-harmonic morphism (see \cite{LO}), the inversion of sphere $\phi:\r^m\setminus\{0\}\longrightarrow \r^m,\;\; \phi(x)=x/|x|^2$ is a biharmonic map if and only if $m=4$ (cf. \cite{BK}), and the radial projection $\phi:\r^m\setminus\{0\}\longrightarrow S^{m-1},\;\; \phi(x)=x/|x|$ is a biharmonic morphism if and only if $m=4$ (cf. \cite{Ou2})!
\item[C3.] Surfaces $M^2\hookrightarrow \r^4$;
\item[C4.] Biharmonic graphs of functions $f:\r^m\supset \Omega\longrightarrow \r$. Notice that Chen's conjecture is local by nature and locally, every hypersurface is given by the graph of a function. For biharmonic graph equation and the equivalent statement of Chen's conjecture for biharmonic graph, see the author's paper \cite{Ou5}.
\end{itemize}
{\bf (2) The generalized Chen's Conjecture (C-M-O '01):} Notice that the Generalized Chen's conjecture for biharmonic submanifolds states that any biharmonic submanifold of nonpositively curved manifold is a minimal one. The counter examples constructed by the author and Tang in \cite{OT} are biharmonic submanifolds which are not complete and the ambient spaces have nonconstant nonpositive sectional curvature. So it would be interesting to know if the following conjectures are true.
\begin{itemize}
\item[GC1.] Conjecture 1: Any biharmonic submanifold of $H^n(-1)$ is a minimal one;
\item[GC2.] Conjecture 2 (\cite{AM}, \cite{Ma6}): Any complete biharmonic submanifold of a nonpositively curved manifold is a minimal one.
\end{itemize}
{\bf (3) Conjectures on biharmonic submanifolds of spheres} (\cite{BMO2}, \cite{BMO}, \cite{BMO3} and \cite{BMO6}):
\begin{itemize}
\item[SC1.] Conjecture 1: The only proper bih. hypersurface of $S^{m+1}$ is a part of $S^m(\frac{1}{\sqrt{2}})$ or $S^p(\frac{1}{\sqrt{2}})\times S^p(\frac{1}{\sqrt{2}})$ with $p+q=m,\;\; p\ne q$;
\item[SC2.] Conjecture 2: Any proper bih, submanifold of $M^m\hookrightarrow S^{m+p}$ has constant mean curvature.
\item[]
\end{itemize}
{\em Final Remarks:} Biharmonic submanifolds ( or more generally, biharmonic maps) has become a popular subject of research in recent years with many significant progresses made by mathematicians around the world. It is impossible (so the author had to give up any attempt) to give a complete account of all results in this brief survey within the page limit of the proceeding. For example, many recent study has been made on biharmonic submanifolds of other spaces more general than the space forms. Interested readers may refer to \cite{IIU1}, \cite{IIU2}, \cite{FLMO} and \cite{Zh} for some classification of biharmonic submanifolds of in $\mathbb{H}P^n$ or $CP^n$; \cite{FO2}
\cite{FO3}, \cite{FO4}, and \cite{Sa} for biharmonic submanifolds in Sasakian space forms; \cite{MU}, \cite{Ur1} and \cite{Ur4} for biharmonic submanifolds in K$\ddot{\rm a}$hler manifolds;  \cite{Ou2}, \cite{OT} and \cite{TO} for biharmonic hypersurfaces in conformally flat spaces; \cite{OSU} \cite{Ou4} for biharmonic hypersurfaces in symmetric or Einstein space; \cite{OW2} for biharmonic surfaces of Thurston's eight 3-dimensional geometries; \cite{FOR} for biharmonic submanifolds with parallel mean curvature vector fields in $S^n\times \r$, and \cite{Ouy}, \cite{Zh2}, \cite{Sa2}, \cite{LD}, \cite{LDZ} and \cite{DO} for biharmonbic submanifolds of pseudo-Riemannian manifolds.

\end{document}